\documentclass[11pt]{article}
\usepackage{latexsym}
\usepackage{amssymb}
\newcommand{\E}{\mathbb E}
\newcommand{\qed}{\hbox{\hskip 6pt\vrule width6pt height6pt depth1pt \hskip1pt}}
\begin{document}

%\title{Emergent Structures in Large Networks}
%\author{David Aristoff$^1$ and Charles Radin$^2$}
\centerline{{\bf Emergent Structures in Large Networks}}

\vskip.2truein
\centerline{David Aristoff$^1$\hskip.7truein Charles Radin$^2$}
\vskip.3truein

\centerline{October, 2011}

\vskip1truein \noindent
Abstract. We consider a large class of exponential random graph models
and prove the existence of a region of parameter space corresponding to
multipartite structure, separated by a phase transition from a
region of disordered graphs.

%\vskip4truein 
\vfill
\phantom{}$^1$ Department of Mathematics

\hskip.1truein University of Minnesota

\hskip .1truein Minneapolis, MN\ \  55455

\phantom{}$^2$ Department of Mathematics

\hskip.1truein University of Texas

\hskip .1truein Austin, TX\ \ 78712
%\date{October 2011}

%\maketitle

\newpage

\medskip 
\centerline{\bf I. Introduction and statement of results} 
\medskip

Complex networks, including the internet, world wide web, social
networks, biological networks etc, are often modeled by probabilistic
ensembles with one or more adjustable parameters; see for instance [N]
and [L],
and the many references therein. We will use one of these standard
families, the exponential random graph models 
(see references in [CD], [F1], [F2], [L] and [N]), to study how
multipartite structure can exist in such networks, stable against
random fluctuations, in imitation of the modeling of crystalline
structure of solids in thermal equilibrium.

We will be considering a large family of exponential random
graph models, but for simplicity we first discuss the particular case
introduced by Strauss [S] in which a graph $G_N$ with $N$ nodes,
$E_N(G_N)$ edges and $T_N(G_N)$ triangles is probabilistically modeled with
the following two-parameter probability mass function:

$$Prob_{\alpha_1,\alpha_2}(G_N)={e^{\alpha_1 E_N(G_N)+\alpha_2 T_N(G_N)}
\over \hbox{normalization}}.\eqno(1)$$

The maximum number of edges in $G_N$ is of order $N^2$, and for
triangles order $N^3$; it will be useful below if we renormalize
quantities. We will work with edge and triangle {\it
densities}, $e_N(G_N)\equiv E_N(G_N)/N^2$ and $t_N(G_N)\equiv T_N(G_N)/N^3$,
and introduce new parameters $\beta_1, \beta_2$ so that 

$$Prob_{\beta_1,\beta_2}(G_N)={e^{N^2[\beta_1 e_N(G_N)+\beta_2 t_N(G_N)]}\over \hbox{normalization}}. \eqno(2)$$

We think of the parameters $\beta_1, \beta_2$ as representing
mechanisms for influencing the network, as pressure and temperature do
in models of materials in thermal equilibrium.
Indeed it is easy to see by
differentiation that
if $\beta_1$ is fixed, varying $\beta_2$ will vary the mean value of
the triangle density; similarly if $\beta_2$ is fixed, varying
$\beta_1$ will vary the mean value of the edge density. Furthermore if
the mean value $\E_{\beta_1,\beta_2}[e_N(G_N)]$ of 
$e_N(G_N)$ is fixed and $\beta_2<< 0$ then, as we
will see below, the random graph will have a very low value for the
mean value $\E_{\beta_1,\beta_2}[t_N(G_N)]$ 
of $t_N(G_N)$. However if 
$\E_{\beta_1,\beta_2}[e_N(G_N)]$
is fixed any variation of $\beta_2> 0$ does not affect 
$\E_{\beta_1,\beta_2}[t_N(G_N)]$ (when $N$ is large) [RY].
It is natural to treat separately the cases $\beta_2<0$
and $\beta_2>0$. The former is called {\it repulsive}, the latter
{\it attractive}; see [RY]. The attractive case $\beta_2>0$ has been completely
analyzed in [RY], so we concentrate here on the case with repulsion, $\beta_2<0$.

It is useful to analyze the phenomenon in the last paragraph, as
regards $\beta_2<<0$,  in two
stages. First, consider the nonprobabilistic optimization problem in
which one maximizes the edge density among those graphs $G_N$ of $N$
nodes which have no triangles, $t_N(G_N)=0$, corresponding intuitively
to $\beta_2=-\infty$. This was solved by Tur\'an
[T], who showed that the optimum is uniquely achieved by the complete
bipartite graph with equal size parts. (The parts differ by 1 if $N$
is odd.) One can understand the Strauss model as a two stage
generalization of this optimization problem. First one considers the 
two-parameter set of graphs:

$${\cal X}_N(e,t)\equiv \{G_N\, :\, e_N(G_N)=e, \
t_N(G_N)=t\}. \eqno(3)$$

\noindent Then one studies the interaction of the two conditions,
$e_N(G_N)=e$ and $t_N(G_N)=t$, through the cardinality 
$|{\cal X}_N(e,t)|$
of 
${\cal X}_N(e,t)$. Specifically, consider the {\it entropy}, defined
on probability mass functions $\rho$ on the set ${\cal X}_N(e,t)$
by:

$$S_N(e,t)[\rho]\equiv -\sum_{j\in
{\cal X}_N(e,t)}\rho_j\ln(\rho_j). \eqno(4)$$

\noindent It is easy to prove that $S_N(e,t)$ is maximized uniquely by
the uniform distribution $\tilde \rho(e,t)$, and that 
$$S_N(e,t)[\tilde
\rho(e,t)]=\ln(|{\cal X}_N(e,t)|). \eqno(5)$$ 

\noindent If one alters
the optimization so that one doesn't restrict $\rho$ to be supported
in ${\cal X}_N(e,t)$ but instead assumes the mean values of the
two densities $e_N(G_N)$ and $t_N(G_N)$ with respect to $\rho$
are fixed, then using
Lagrange multipliers $N^2\beta_1$ for $e_N(G_N)$ and $N^2\beta_2$ for
$t_N(G_N)$ in the optimization of the entropy leads to the unique optimizer given in (2), in which the
$\beta$'s control the mean values. It was shown in [CD] that if $\beta_2<<0$ then 
in a technical sense $G_N$ for large $N$ looks like the complete
bipartite graph with equal parts, with some edges
randomly removed, but in particular $t_N(G_N)\approx 0$. On the
other hand we will see that when $-1/3<\beta_2<0$, the edges in $G_N$ are
roughly independent, so fixing $e_N(G_N)$ automatically fixes
$t_N(G_N)$, not leaving any flexibility.

Although the networks corresponding to $\beta_2>0$ do not have
interesting structure, there is still an interesting phenomenon in
this regime,
associated with sensitivity to variation of the parameters. More
specifically, it was proven in [CD] that at certain values of
$\beta_2>0$ there is a special value $\beta_1=\beta_1(\beta_2)$ at which small
changes of $\beta_1$ with $\beta_2$ fixed
lead to a jump in the mean value of the density
$e_N(G_N)$. Furthermore, it was shown in [RY] that {\it all} singular
behavior of the distributions is concentrated on a certain curve $\beta_1=q(\beta_2)$. 
(We will clarify the meaning of ``singular'' below.)

We are interested here in the more complicated case $\beta_2<0$. 
As noted above, if $|\beta_2| < 1/3$ then for large $N$ the graph 
will have approximately independent edges; in particular we will show
below that the difference

$$ \Big(\E_{\beta_1,\beta_2}[e_N(G_N)]\Big)^3-\E_{\beta_1,\beta_2}[t_N(G_N)]
\eqno(6)$$ 

\noindent has limit 0 as $N\to \infty$. However one might expect
from Tur\'an's theorem 
that for any fixed $\beta_1$ (or mean value of
$e_N(G_N)$), once $\beta_2$ is sufficiently negative the graph
should look bipartite, and so the difference should be roughly
$\displaystyle (\E_{\beta_1,\beta_2}[e_N(G_N)])^3\ne 0$. We prove
this below but furthermore prove that the qualitative structural 
change occurs abruptly:
in order to accomplish the change, for each $\beta_1$ the distribution
exhibits ``singular behavior'' at some $\beta_2 <0$. Before clarifying
the meaning of ``singular'' we generalize the role played by triangles
in the Strauss model to the following two-parameter {\it exponential
  random graph model}:

$${\Bbb P}_{\beta_1,\beta_2}(G_N)=
e^{N^2[\beta_1 t(H_1,G_N)+\beta_2
  t(H_2,G_N) - \psi_N(\beta_1,\beta_2)]}, \eqno(7)$$

 \noindent where: $H_1$ is an edge,
 $H_2$ is any finite simple graph with $k\ge
  2$ edges, and $t(H_i,G_N)$ is
  the density of graph homomorphisms $H \to G_N$:
  
 $$ t(H_i,G_N)={|\hbox {hom}(H_i,G_N)|\over |V(G_N)|^{|V(H_i)|}}, \eqno(8)$$

\noindent with $V(\cdot)$ denoting a vertex set. The term $\psi_N(\beta_1,\beta_2)$ 
gives the probability normalization. 

Fundamental to our
results are questions of analyticity of the normalization in (7),
which we discuss next. (See [KP] for elementary properties of real
analytic functions of several real variables.) 
An explicit formulation of the normalization is:

$$\psi_N(\beta_1,\beta_2)\equiv 
{1\over N^2}\ln\Big(\sum_{G_N} 
e^{N^2[\beta_1 t(H_1,G_N)+\beta_2 t(H_2,G_N)]}\Big). \eqno(9)$$

\noindent It is proven in [CD] that 

$$\psi_\infty(\beta_1,\beta_2)\equiv \lim_{N\to\infty}\psi_N(\beta_1,\beta_2) \eqno(10)$$

\noindent exists for all $\beta_1,\beta_2$. It is also noted in [RY]
that at points where $\psi_\infty$ is analytic,

$$\frac{\partial}{\partial \beta_j}\psi_\infty(\beta_1,\beta_2)= \lim_{N\to\infty}\frac{\partial}{\partial \beta_j}\psi_N(\beta_1,\beta_2), \eqno(11)$$

\noindent that is, the partial derivatives commute with the limit $N\to \infty$. 
Partial derivatives of $\psi_\infty$, when they exist, 
give information on
the large-$N$ mean and variance of the densities $t(H_1,G_N)$ and $t(H_2,G_N)$ (see [RY])
and it is standard in the corresponding modeling of materials, in part for this reason, to define
phases and phase transitions as follows (see [FR]).

\medskip
 \noindent {\bf Definition.} A phase is an open connected region 
of the parameter space
  $\{(\beta_1,\beta_2)\}$ which is maximal for the condition that
  $\psi_{\infty}(\beta_1, \beta_2)$
  is analytic. 
  There is a phase transition at
   $(\beta_1^\ast, \beta_2^\ast)$ if $(\beta_1^\ast,
   \beta_2^\ast)$ is a boundary point of an
  open set on which $\psi_{\infty}$ is analytic, but $\psi_{\infty}$ is not analytic at 
$(\beta_1^\ast, \beta_2^\ast)$.
\medskip 

So in the above, ``singular'' meant nonanalytic. In this notation our main result is: 
\medskip

 \noindent {\bf Theorem.} Assume the chromatic number $\chi(H_2)$ of $H_2$ is at
least 3. Then there is a curve $\beta_2=s(\beta_1),\ -\infty < \beta_1<
\infty$, in the lower half plane ($\beta_2<0$),
such that the model exhibits a phase transition on the curve. 

\bigskip

\centerline{{\bf II. Proof of the theorem}}
\medskip

 Let $k$ be the number of edges in $H_2$. 
We write $\Bbb P$ for the probability mass function
${\Bbb P}_{\beta_1,\beta_2}$ given by equation (7), and $\Bbb E$ for
the expectation ${\Bbb E}_{\beta_1,\beta_2}$.
\medskip

By Theorem 6.1 in [CD], the analyticity method of the proof of 
Theorem 3.10 in [RY] can be immediately extended to prove 
that $\psi_\infty(\beta_1,\beta_2)$ 
is analytic in the real variables $\beta_1$ and $\beta_2$ when 
$|\beta_2| < 2/[k(k-1)]$. Our proof will be by contradiction, 
so we assume from here on that $\psi_\infty(\beta_1,\beta_2)$ is 
analytic in $\beta_1$ and $\beta_2$ 
on the {\it entire} half line $L = \{(\beta_1^*, \beta_2): \beta_2 < 0\}$, 
where $\beta_1^*$ is arbitrary but fixed. We will find a contradiction, which
will prove the existence of the curve $\beta_2=s(\beta_1)$.

Consider the 
function 
$$C(\beta_1,\beta_2) := 
\left(\frac{\partial \psi_\infty}{\partial
    \beta_1}(\beta_1,\beta_2)\right)^k 
- \frac{\partial \psi_\infty}{\partial
      \beta_2}(\beta_1,\beta_2)
\eqno(12)$$

\noindent Note that $C(\beta_1,\beta_2)$ is analytic on $L$, since $\psi_\infty(\beta_1,\beta_2)$ is.

Proposition 3.2 in [RY] proves, for all $\beta_2 < 0$,
there is a unique solution $u^*(\beta_1,\beta_2)$ to the
optimization of
$$\beta_1 u + \beta_2 u^k - \frac{1}{2} u \ln u -
\frac{1}{2}(1-u)\ln(1-u) \eqno(13)$$

\noindent for $u \in [0,1]$. Then from Theorems
6.1 and 4.2 in [CD] we can use the same argument as used to prove
equations (33) and (34) in [RY] to prove, for $-2/[k(k-1)]<\beta_2<0$:

$$\frac{\partial}{\partial
  \beta_1}\psi_{\infty}(\beta_1, \beta_2)=
\lim_{N\to \infty}{\Bbb E}\{t(H_1,G_N)\}=
  t(H_1,u^*)=u^*(\beta_1, \beta_2), \eqno(14)$$

$$\frac{\partial}{\partial
  \beta_2}\psi_{\infty}(\beta_1, \beta_2)=
\lim_{N\to \infty}{\Bbb E}\{t(H_2,G_N)\}=
  t(H_2,u^*)=\left(u^*(\beta_1, \beta_2)\right)^k. \eqno(15)$$

\noindent It follows that $C(\beta_1^*, \beta_2) =
\left(t(H_1,u^*)\right)^k-t(H_2,u^*)\equiv 0$
for $|\beta_2| < 2/[k(k-1)]$. Since a function of one variable which is
analytic on $L$
and constant on a subinterval must be constant on $L$, it follows that

$$C(\beta_1^*, \beta_2) \equiv 0 \hbox{ on } L. \eqno(16)$$

We next obtain a contradiction to (16), but first we need some
notation; see [CD], [BCL] and [L] for discussions of the ideas behind these
terms, which basically provide the framework for ``infinite volume limits'' 
for graphs, in analogy with the infinite volume limit in statistical mechanics [R].

To each graph $G$ on $N$ nodes we associate the following function on
$[0,1]^2$:

$$f^G(x,y)= 1 \hbox{ if } (\lceil Nx\rceil, \lceil Ny\rceil) \hbox{ is an edge of
  } G, \hbox{ and }
f^G(x,y)=0 \hbox{ otherwise}. \eqno(17)$$

\noindent We define ${\cal W}$ to
be the space of measurable functions $h:[0,1]^2 \to [0,1]$ which
are symmetric: $h(x,y)=h(y,x)$, for all $x,y$. 
For $h\in {\cal W}$ we define

 $$t(H,h)\equiv \int_{[0,1]^\ell} \prod_{(i,j)\in E(H)}h(x_i,x_j)\,dx_1\cdots
   dx_\ell. \eqno(18) $$

\noindent where $E(H)$ is the edge set of $H$, and $\ell = |V(H)|$ is the number 
of nodes in $H$, and note that for a graph $G$, $t(H,G)$ defined in
(8) has the same value as $t(H,f^G)$.
We define an equivalence relation on ${\cal W}$ as follows: $f\sim g$
if and only if $t(H,f)=t(H,g)$ for every simple graph $H$.
Elements of the quotient space, $\tilde {\cal W}$, are called
``graphons'', and the class containing $h\in {\cal W}$ is denoted
$\tilde h$.
The space $\tilde W$ is compact [L].

On $\tilde {\cal W}$ we define a metric in steps as follows.
First, on ${\cal W}$ we define

$$d_\square(f,g)\equiv \sup_{S,T\subseteq [0,1]}\Big| \int_{S\times
  T}[f(x,y)-g(x,y)]\, dxdy\Big|. \eqno(19)$$

\noindent Let $\Sigma$ be the space of measure preserving bijections $\sigma$ of
$[0,1]$, and for $f$ in ${\cal W}$ and $\sigma\in \Sigma$ define
$f_\sigma(x,y)\equiv f(\sigma(x),\sigma(y))$. Using this we define
a metric on $\tilde W$ by

$$\delta_\square(\tilde f,\tilde g)\equiv \inf_{\sigma_1,\sigma_2}
d_\square(f_{\sigma_1},g_{\sigma_2}). \eqno(20)$$

Next we need a few terms associated with $\psi_\infty$. Define on $[0,1]$:

$$I(u)\equiv \frac{1}{2}u\ln(u) + \frac{1}{2}(1-u)\ln(1-u) \eqno(21)$$
and on $\tilde {\cal W}$:

$$I(\tilde h)\equiv \int_{[0,1]^2}  I(h(x,y))\, dxdy. \eqno(22)$$

\noindent Also on $\tilde {\cal W}$ we define:

$$T(\tilde h)\equiv \beta_1 t(H_1,h)+\beta_2 t(H_2,h). \eqno(23)$$

The above are relevant because it is proven in Theorem 3.1 of [CD]
that $\psi_\infty(\beta_1,\beta_2)$ is the solution of an optimization
problem:

$$\psi_\infty(\beta_1,\beta_2)=\sup_{\tilde h\in \tilde {\cal W}}
\big[T(\tilde h)-I(\tilde h)\big]. \eqno(24)$$

\noindent Furthermore, from Theorem 3.2 of [CD] one has some control on the
asymptotic behavior as $N\to \infty$:

$$\delta_\square[\tilde G_N,\tilde F^\ast(\beta_1,\beta_2)]\to 0
\hbox{ in probability as }N\to \infty, \eqno(25)$$

\noindent where $\tilde F^\ast(\beta_1,\beta_2)$ is the (compact) subset of 
$\tilde {\cal W}$ on which $T-I$ is maximized, and $\tilde G_N\equiv
\tilde f^{G_N}$.

We now return to our proof. Fix $\epsilon>0$ and $i\in \{1,2\}$. Recall $\beta_1 = \beta_1^*$ is 
fixed arbitrarily. Write $\tilde{F}^*(\beta_2)$ for the set $\tilde{F}^*(\beta_1,\beta_2) \subset \tilde W$ defined above.
Using Theorem 7.1 in [CD], 
choose $\beta_2^\prime$ sufficiently negative so that for every $\beta_2<\beta_2^\prime$
$$\sup_{\tilde{f} \in \tilde{F}^*(\beta_2)} \delta_\square(\tilde{f},p\tilde{g}) < \frac{\epsilon}{3k}. \eqno(26)$$
Using Theorem 3.2 in [CD], choose $N_0(\beta_2)$ such that $N > N_0(\beta_2)$ implies 
$${\Bbb P}\left[\delta_\square(\tilde{G}_N,\tilde{F}^*(\beta_2))\ge \frac{\epsilon}{3k}\right]<\frac{\epsilon}{3k}. \eqno(27)$$

\noindent
Let $A_{\epsilon,N} = \{G_N: \delta_\square (\tilde G_N, \tilde{F}^*(\beta_2)) < \epsilon/(3k)\}$. 
By compactness of $\tilde{F}^*(\beta_2)$ we may 
choose $\tilde{h}_{G_N} \in \tilde{F}^*(\beta_2)$ corresponding to each $G_N \in A_{\epsilon,N}$ such  
that $$\delta_\square(\tilde{G}_N,\tilde{h}_{G_N}) < \frac{\epsilon}{3k}.\eqno(28)$$

\noindent
Write ${\Bbb E}\big|_A$ for the 
restriction of the expectation to the set $A$. Using (26) and (28) we have that   
$${\Bbb E}\big|_{A_{\epsilon,N}}[\delta_\square(\tilde{G}_N,p\tilde{g})] 
= \sum_{G_N \in A_{\epsilon,N}}\delta_\square(\tilde{G}_N,p\tilde{g}){\Bbb P}(G_N)$$
$$ \le \sum_{G_N \in A_{\epsilon,N}}\left[\delta_\square(\tilde{G}_N,\tilde{h}_{G_N}) + \delta_\square(\tilde{h}_{G_N},p\tilde{g})\right]{\Bbb P}(G_N)$$
$$ < \sum_{G_N \in A_{\epsilon,N}} \left[\frac{\epsilon}{3k} + \frac{\epsilon}{3k}\right]{\Bbb P}(G_N) \le \frac{2\epsilon}{3k} \eqno(29)$$

\noindent
for $N > N_0(\beta_2)$. Now write $\bar{A}_{\epsilon,N}$ for the complement of $A_{\epsilon,N}$. 
Then by Lemma 3.12 in [RY], equation (29), and the fact that $\delta_\square(\cdot,\cdot) \le 1$, 
$$\Big|{\Bbb E}\left[t(H_i,G_N)\right] - t(H_i,pg)\Big|$$
$$\le{\Bbb E}\left[\left|t(H_i,G_N)-t(H_i,pg)\right|\right]$$
$$ \le k \cdot {\Bbb E}\left[\delta_\square(\tilde{G}_N,p\tilde{g})\right]$$
$$ = k\cdot\left({\Bbb E}\big|_{A_{\epsilon,n}}\left[\delta_\square(\tilde{G}_N,p\tilde{g})\right]+ {\Bbb E}\big|_{\bar{A}_{\epsilon,N}}\left[\delta_\square(\tilde{G}_N,p\tilde{g})\right]\right)$$
$$ \le k\cdot\left(\frac{2\epsilon}{3k}+\frac{\epsilon}{3k}\right) = \epsilon \eqno(30)$$

\noindent
for $N>N_0(\beta_2)$. Using the identity 
$$\frac{\partial \psi_N}{\partial \beta_i}(\beta_1^*,\beta_2) = {\Bbb E}\left[t(H_i,G_N)\right] \eqno(31)$$
along with (11), we 
may take the limit $N\to \infty$ in (30) to obtain 
$$\left|t(H_i,pg)-\frac{\partial \psi_\infty}{\partial \beta_i}(\beta_1^*,\beta_2)\right| < \epsilon. \eqno(32)$$
Since $\epsilon>0$ was arbitrary, 
$$\lim_{\beta_2\to -\infty} \frac{\partial \psi_\infty}{\partial \beta_i}(\beta_1^*,\beta_2) = t(H_i,pg). \eqno(33)$$

\noindent
Direct computation using equation (2.10) in [CD] yields:
$$t(H_2,p g)=0\ \hbox{ and }\ t(H_1,p
g)=\frac{e^{2\beta_1}(\chi(H)-2)}{(1+e^{2\beta_1})(\chi(H)-1)}>0. \eqno(34)$$
Now, by combining (12) with (33)-(34) we find
$\lim_{\beta_2\to -\infty} C(\beta_1^*, \beta_2) >0$, in contradiction with (16), which
proves the theorem. \qed
\vskip.2truein

\centerline{{\bf III. Conclusion}} 
\medskip

Consider any of the two-parameter exponential random
graph models with repulsion covered by our theorem. Define the `high
energy phase' of the parameter space $\{(\beta_1,\beta_2)\,| \,
\beta_2<0\}$ as that domain of analyticity of
$\psi_\infty(\beta_1,\beta_2)$ which contains the strip
$-2/[k(k-1)]<\beta_2< 0$. The order parameter $C(\beta_1,\beta_2)$ is
identically zero in this phase (as one can see for instance by
connecting by an analytic curve any given point in the open, connected
phase to a point in the strip, and complexifying). We have proven that this
phase is separated from the low energy regime in the sense that for
each $\beta_1$ there is some $\beta_2^\prime$ such that the segment
$\{(\beta_1,\beta_2)\,|\, \beta_2< \beta_2^\prime\}$ does not
intersect the phase. Our proof is based
on the traditional modeling of equilibrium statistical mechanics using
analyticity and an order parameter
[R], [K], [Y]. And we emphasize that this method
could not have been used to prove the transition found in [RY] for attractive 
exponential random
graph models since there is a critical point for that
transition: indeed there is only one phase.

In comparison with traditional models from statistical mechanics,
exponential random graph models could be thought of as either infinite
range, or infinite dimensional, which suggests a relation with `mean
field theories' [K], [Y]. Mean field theories grew out of the work of
van der Waals, who obtained a general description of fluids by adroitly
replacing the interaction between each molecule and the rest of the
fluid by an average or mean field, among other things losing track of
the spatial separation of the interacting molecules. This proved to be
a useful approximation to understand gas/liquid phase transitions in
which the most relevant part of the particle interaction is a (long
range) attraction. It is not too surprising that exponential random
graph models with attractive interaction could therefore yield a phase
transition like that of the liquid/gas transition, as was shown in
[RY]. In the present paper we obtain a transition more like a
fluid/solid transition, in which there is a change of `symmetry' from
disordered to multipartite. It is less intuitive to use a long range
repulsion to model a solid/fluid transition, so the materials analogy
of our models with repulsion is less compelling than for our models
with attraction. Therefore the relation of these models with repulsion to
mean field approximations might be particularly illuminating. 

There remain many open questions. Perhaps the most pressing is the
character of the singularity of $\psi_\infty(\beta_1,\beta_2)$ at the
boundary of the high energy phase. In the attractive case there is
only one phase but there are jump discontinuities, in the first
derivatives of $\psi_\infty(\beta_1,\beta_2)$ (namely the average edge
and energy densities), across a curve where two regions of the phase
abut, while the edges are independent throughout the phase [RY]. We do
not know the nature of the singularity at the boundary of the high
energy phase for the case of repulsion studied in this paper, though
we expect the first derivatives of $\psi_\infty(\beta_1,\beta_2)$ to
be discontinuous across the boundary. In analogy with equilibrium
materials there may be multipartite phases with different numbers of
parts at low energy, though this may require more complicated
interactions [CD].

\bigskip

\noindent Acknowledgments. It is a pleasure for CR to acknowledge useful
discussions with Mei Yin, and support at a workshop of The American
Institute of Mathematics in August 2011.
\bigskip

\bigskip

\centerline{Bibliography}
\bigskip
\noindent {[BCL]} J.T. Chayes, C. Borgs and L. Lovasz,
Moments of two-variable functions and the uniquenes of graph limits,
GAFA 19 (2010) 1597-1619.
\medskip

\noindent {[CD]} S. Chatterjee and P. Diaconis, Estimating and understanding
exponential random graph models, arxiv:1102.2650v3
\medskip

\noindent {[F1]} S.E. Fienberg, Introduction to papers on the
  modeling and analysis of network data, {Ann. Appl. Statist.} 4 (2010) 1-4.
\medskip

\noindent {[F2]} S.E. Fienberg, Introduction to papers on the
  modeling and analysis of network data II, {Ann. Appl. Statist.}
  4 (2010) 533-534.
\medskip

\noindent {[FR]} M.E. Fisher and C. Radin, Definitions of thermodynamic phases and
phase transitions, 2006 workshop report,
http://www.aimath.org/WWN/phasetransition/Defs16.pdf
\medskip

\noindent {[K]} L.P. Kadanoff, Theories of matter: infinities and renormalization,
The Oxford Handbook of the Philosophy of Physics, ed. Robert
Batterman, Oxford University Press, 2011.
\medskip

\noindent [KP] S.G. Krantz and H.R. Parks, A Primer of Real Analytic Functions,
$2^{nd}$ ed., Birkhauser, 2002.
\medskip

%\noindent [LS] L. Lov\'{a}sz, and B. Szegedy, Limits of
%dense graph sequences, {J. Combin. Theory Ser. B} 98 (2006)
%933-957.

\noindent [L] L. Lov\'{a}sz, Very large graphs, Current Developments
in Mathematics, 67-128 (2008).

\noindent {[N]} M.E.J. Newman, Networks: an Introduction, Oxford University Press,
2010.
\medskip

\noindent {[R]} D. Ruelle, { Statistical Mechanics; Rigorous Results}, Benjamin, New
York, 1969.
\medskip

\noindent {[RY]} C. Radin and M. Yin, Phase transitions in exponential random
graphs, arxiv:1108.0649
\medskip

\noindent {[S]} D. Strauss, On a general class of models for interaction, SIAM
Rev. 28 (1986) 513-527.
\medskip

\noindent {[T]} P. Tur\'an, On an extremal problem in graph theory (in
Hungarian), Matematikai \'es Fizikai Lapok 48 (1941) 436–452.
\medskip

\noindent {[Y]} J.M. Yeoman, Statistical Mechanics of Phase Transitions, Clarendon
Press, Oxford, 1992.

\end{document}